\documentclass[12pt]{article} 
\usepackage{amssymb,epsfig,euscript,latexsym, amsmath}

\makeatletter
\@addtoreset{figure}{section}
\def\thefigure{\thesection.\@arabic\c@figure}
\@addtoreset{equation}{section}

\makeatother

\def\fder#1#2{{\frac {d #1}  {d #2}}}
\def\pder#1#2{{\frac {\partial #1} {\partial #2}}}

\def\eqref#1{(\ref{#1})}

\def\ex{\exp \frac{1}{\hbar}}          

\def\ctablec#1{\vbox {\offinterlineskip \hrule\halign
{\vrule\enspace\hfil $##$\hfil\strut\enspace\vrule&&\enspace\hfil
$##$\hfil\enspace\vrule\cr#1}\hrule}}

\begin{document}

\title{A Rigorous Path Integral Construction in Any Dimension}

\author{Alexander Dynin} 

\date{ Department of Mathematics, Ohio State University \newline  
Columbus, OH 43210}

\maketitle 

\begin{abstract}
A rigorous Path Integral construction 
for a wide class of Weyl evolution operators
 is based on a
pseudo-differential $\Omega $-calculus on  flat phase spaces of finite
and infinite dimensions. 
\end{abstract}

\section{Introduction}
Among the various uses of the Path Integrals the earliest and foremost is 
the application to  quantum evolution. 
The Path Integral formalism is appreciated
 because of its compact operational notation 
with all the gamut  of the finite dimensional Integral Calculus, such 
as integration by parts,  repeated integration, 
canonical substitutions, analytic continuation (Wick rotation),
stationary phase approximations etc.

 However,
in R.Feynman words, the Path Integral is ``an intuitive leap at 
mathematical
formalism''. 
A natural  justification  would be by  suitable
integral approximations, the route chosen originally by Feynman himself 
in the 40's via time-slicing and discretization. Unfortunately the 
discretization 
ambiguities along with the convergence problems have plagued  the deed 
from the start
(the notable exception~\cite{Nel} is for  rather special hamiltonians).
 
 \smallskip  

Nevertheless we propose a \emph{rigorous time-slicing}
construction of the (flat) phase space 
Path Integral for propagators both in Quantum Mechanics and Quantum Field 
Theory 
  for a fairly general class  of quasi-dissipative  quantum observables  
  (e.g. the Schr\H{o}dinger hamiltonians with smooth scalar potentials of 
any 
  power growth). Moreover we allow  
  time-dependent hamiltonians and  a great variety of 
  discretizations, in particular,  the standard, Weyl, and normal ones.

\paragraph{Abstract Cauchy Problem.} Consider the Initial Value Problem
$$\fder{\psi}{t} + A(t)\psi (t) = 0,\quad \psi (0) = \psi_0,\quad 
0\leq t\leq T, \quad \psi (t) \in \mathcal H ,$$
wherein $\mathcal H$ is a Hilbert space and  
$A(t)$ is a family of (usually)  unbounded operators on $\mathcal H$.

The Cauchy Problem is called {\it proper} relative a dense subspace 
$\mathcal S$ in $\mathcal H$
if there is the unique solution $\psi$ for every $\psi_0 \in \mathcal S$ 
and the  {\it Evolution Operator }
$$U(t'',t')\psi (t')=\psi (t''),\quad \psi (t') \in {\mathcal S}, 
\quad t''>t',$$
is bounded on $\mathcal H$.

The Evolution Operator may be sought in the time-slicing form of the 
 strong operator limit \emph{Product Integral}:
\begin{align*}
U(t'',t') &= \prod_{t''\ge t \ge t'}\exp [-A(t)dt] \\
          &:= \lim_{|{\mathcal P}|\rightarrow  0}U^{\mathcal P}(t'',t')
: =\lim_{|{\mathcal P}|\rightarrow 0} \prod_{t_{j+1}>t_j}
 \exp [-A(t_j) 
\Delta t_j].
\end{align*}
Here $\mathcal P$ is a finite partition 
$0 \leq t'=t_0 < \dots <t_j<t_{j+1}< \dots <t''=t_ \mathcal P \leq T$ of 
the interval 
$t' \leq t \leq t'', \   \Delta t_j=t_{j+1}-t_j,\     
| \mathcal P|= \max_{j} | \Delta t_{j}|$ 

\paragraph{From Product to Path Integrals on the Phase Space.}

Consider a quantum evolution equation on $\mathcal L^2 (\mathbb R^d)$
$$
\pder{\psi}{t}(t,q)+\frac{i}{\hbar}f(t,q,\frac{i}{\hbar}\pder{}{q})\psi 
(t,q)=0,\quad \psi (0,q)=\psi_0,
$$
with a pseudodifferential operator $f(t,q,\frac{i}{\hbar} 
\pder{}{q})$ in the standard form of the $qp-$quantization of a 
complex-valued 
function $f(t,q,p)$ on the phase space $\mathbb R^{2d}$, the standard 
symbol of $f(t,q,\frac{i}{\hbar} \pder{}{q})$. \smallskip

 The Tobocman's version of the Dirac-Feynman  Ansatz:   
{\it For small $\Delta t $ the standard symbol $<p|U(t+\Delta t,t)|q>$ of 
the propagator \mbox{$U(t+\Delta t,t)$} is approximately equal to 
$\exp [-\frac{i}{\hbar} f(t,q,p)\Delta t]$.}  \smallskip
 Then according to the product rule for the standard symbols, the 
standard symbol of $U^{\mathcal P}$ is approximately equal to the 
distributional multiple integral
$$
\int \prod_{j=1}^{d-1} d\lambda _{\hbar} (t_j) \exp 
\frac{i}{\hbar}\sum_{j=0}^{d-1}[ p(t_{j+1})\Delta q(t_j)-f(t_j,q(t_j), 
p(t_{j+1})\Delta t_j]
$$
where 
$d\lambda _{\hbar} (t_j)=(2\pi \hbar)^{-d} dq(t_j)dp (t_j)$ 
is the Lebesgue-Liouville measure on the phase space ${\mathbb R}^{2d}$. 
\smallskip 

As the mesh $|\mathcal P |\rightarrow 0$, the multiple integrals are 
presumed to converge to a Hamiltonian
Path Integral for the standard symbol of the evolution operator $U(t'',t')$
 $$
\int \prod_{t''\geq t\geq t'} d\lambda _{\hbar} (t) \exp \frac{i}{\hbar} 
\int_{t'}^{t''} [p(t)
\dot{q} (t)-f(t,q(t),p (t)]dt
$$
where 
$d\lambda _{\hbar} (t) = 
(2\pi \hbar)^{-d} dq(t)dp (t)$
is  {\it ''the Feynman-Liouville measure''} on the  space of
paths from $q(t')=q$ to $p(t'') = p $ in $\mathbb R^{2d}$ and 
$$\int_{t'}^{t''} [p(t)\dot{q} (t)-f(t,q(t),p(t)]dt$$ 
is the hamiltonian symplectic action functional on that space.\smallskip

The standing physical presumption:  all  calculus rules are valid in the 
limit.
 However there are two fundamental mathematical problems: \newline 
 \emph{The 
 validity of the DFT-ansatz} and \emph{the existence of the limit.} So far
  both problems have been settled only for  special  $f$
  (c.f.~\cite{Nel,Mor,Hida}).

\paragraph{Euler Detour.}

Euler polygonal approximation for the solution of the  quantum Cauchy 
Problem
$$
\fder{\psi}{t}+\hat {f}(t)\psi (t)=0
$$
is the finite difference approximation  
$$
\frac{\psi (t_{j+1})- \psi(t_j)}{\Delta t_j}+
\hat {f}(t_j)\psi (t_j)=0, \quad \psi (t_0)=\psi _0, 
$$
or
$$\psi (t_{j+1})=( 1 - \hat {f}(t_j)\Delta t_j)\psi (t_j),
$$
so that the Evolution Operator might be  the strong operator limit
$$
U(t'',t')=\lim_ {|\mathcal P|\rightarrow 0}
 \prod_j ( 1 - \hat {f}(t_j)\Delta t_j) 
         :=\prod_{t''\geq t\geq t'}[ 1 - \hat {f}(t)dt].
$$

If  $\hat {f}(t)$ is a standard pseudodifferential operator then 
the partial product approximations are pseudodifferential operators again. 
However, if the order of $\hat {f}(t)$ is positive, then the order of 
these 
pseudodifferential operator 
approximations increases to infinity and the convergence of their symbols 
is out of the  control.\smallskip

Fortunately, the backward Euler approximation
\[
\frac{\psi(t_{j+1})-\psi (t_j)}{\Delta t_j} +
 A(t_{j+1})\psi (t_{j+1})=0
\] 
suggests
$$
U(t'',t')=\prod_{t''\geq t\geq t'}
( 1 +\hat {f}(t)dt)^{-1}
$$
with zero order approximation symbols.\smallskip

 Our main result is  that for \emph{apt} functions $f$ this backward 
approximation
 entails (in 
the spirit of the  DFT ansatz) 
\[
U(t'',t')=\prod_{t'\geq t\geq t'} \left [( 1 - f(t)dt)^{-1} 
\right] \widehat{}
\]
leading to a Path Integral representation of the symbol $<\! 
p|U(t'',t')|q\! >$. 
 \smallskip 
  
  Incidentally, the Green function (the coordinate propagator) can be 
easily 
  expressed via the symbol:
  $$<\! q''|U(t'',t')|q'\! >=\int dp<\!q''|p\! ><\! p|U(t'',t')|q'\! >.$$

\section{ Rigorized $\Omega$-symbolic calculus}
This section provides necessary technical tools.

\smallskip

  For $z=(q,p) \in \mathbb R^d \times \mathbb R^d$ introduce the complex 
coordinates
   $z^+=2^{-1/2}(q+ip), \quad z^-=2^{-1/2}(q-ip)$
 so that  the standard 
symplectic form $[(p_1,q_1), (p_2, q_2)]=p_1q_2-p_2q_1 $ on   $\mathbb 
R^{2d}$
becomes
\[
 \frac{1}{i} [z_1,z_2]=
\frac{1}{i} (z_{1}^{+}z_{2}^{-}
-z_{1}^{-} z_{2}^+).
\]

The \emph{$\hbar$ -Symplectic  Fourier  transform}
is defined as
 \[ {\tilde{f}}(\zeta)=\int f(z)
e^{\frac{1}{\hbar}[z,\zeta]}d\lambda _{\hbar}(z),\quad  
 d_{\hbar} \lambda (z)  = \frac{1}{(\pi \hbar )^d } dz^+ dz^-.
 \]

\paragraph{Heisenberg Canonical Commutation Relations} 
$z \mapsto \hat z$ over $\mathbb R^{2d}$ in a 
Hilbert space $\mathcal  H$ is a linear map of $\mathbb R^{2d}$ to 
essentially self-adjoint operators on a common invariant subspace 
$\mathcal G$ of  
$\mathcal  H$ (the G\H{a}rding domain) such that 
$[\hat{z} _1,{\hat{z} _2}]=\hbar [z_1,z_2] \mathbf 1$. \smallskip 

E.g., for  the Schr\H{o}dinger (position) representation
on ${\mathcal L}^2(\mathbb R^d)$ 
\[ 
 \hat{z} \psi (x)=(qx)\psi (x) + \frac{\hbar}{i} \pder{\psi }{x} (x),
\]
the Schwartz space $\mathcal S (\mathbb R^d)$  may be chosen as
a G\H{a}rding domain. \smallskip 

Other examples are the momentum or mixed momentum-position representations,
holomorphic Bargmann-Segal representation (conducive to the coherent 
states 
Path Integral), Gelfand-Zak representation
in the Solid State Physics, the Cartier compact representation
in $\theta $-functions. \smallskip

By  the von Neumann-Stone theorem, for given $\hbar >0$ 
any Heisenberg Canonical Commutation Relations is unitary equivalent 
to a direct sum of the Schr\H{o}dinger representations. Thus we may chose 
the G\H{a}rding domain $\mathcal G (\mathcal H)$ to be unitary equivalent 
to a direct sum of the spaces $\mathcal S (\mathbb R^d)$.

 Correspondingly,  the dual space $\mathcal G' (\mathcal H)$ is unitary 
 equivalent to a direct sum of the spaces $\mathcal S' (\mathbb  R^d)$.
  
\paragraph{Weyl operators $\hat{f}$ on $\mathcal H$} 
associated with  generalized functions $f\in  \mathbb R^{2d}$ are 
continuous linear operators from $\mathcal G (\mathcal H)$ to $\mathcal G' 
(\mathcal H)$
\[
\hat{f} = \int \tilde{f} (\zeta ) \ex [\zeta,\hat{z} ] d\lambda _{\hbar}  
(\zeta )
\]
wherein
\[
[\zeta ,\hat{z}] = \zeta ^+ \hat{z}^- - \zeta ^-\hat{z}^+.
\]
A version of the {\it Schwartz Kernel Theorem} states that a linear 
operator from $\mathcal G (\mathcal H)$ to $\mathcal G' (\mathcal H)$
is continuous if and only if it is a Weyl operator $\hat{f}$.

\paragraph{$\Omega$ -symbols.}
Consider a formal power series over $\mathbb C$ 
$$
\Omega (\zeta ) = 1 + \sum_{|\alpha | >0} c_\alpha z^\alpha.
$$
A  {\it formal} $\Omega $-{\it symbol} of $f\in {\mathcal S'}(\mathbb 
R^{2d})$ is the 
formal power series over $\mathcal S'(\mathbb R^{2d})$ 
defined via
\[
\tilde{f} ^{\Omega} (\zeta ) = \tilde{f}/\Omega (\zeta).
\]
Obviously, this has sense  for polynomial $f(z)$ when various 
$\Omega$  provide common
ordering rules according to the following table (c.f.~\cite{Aga}):

 \[                   
\ctablec{\textrm{Name} &\Omega (\zeta ) &\textrm{Ordering}~(d=1)\cr 
\noalign{\hrule}
&&\cr
\textrm{Weyl} &1 &q^np^m\leftrightarrow\cr &\textrm{} &\frac{1}{2^n}
\sum_{j=0}^{n} \left (\! \begin{array}{c} n\\j \end{array} \! \right ) 
\hat{q}^{n-j} \hat{p}^m\hat{q}^j\cr
&&\cr
\noalign{\hrule}
&&\cr
\textrm{Standard} &e^{\frac{1}{4}[(\zeta ^+)^2-(\zeta ^-)^2)} 
&q^np^m\leftrightarrow 
\hat{q}^n\hat{p}^m\cr
(qp~\textrm{or~Kohn-Nirenberg}) &\textrm{} &\cr &&\cr
\noalign{\hrule}
&&\cr
\textrm{Antistandard} &e^{-\frac{1}{4}[(\zeta ^+)^2-(\zeta ^-)^2)]} &q^np^m
\leftrightarrow \hat{p}^m\hat{q}^n\cr
(\textrm{or}~pq) &\textrm{} &\cr &&\cr
\noalign{\hrule}
&&\cr
\textrm{Normal} &e^{\frac{1}{2}\zeta ^+ \zeta ^-} &(z^+)^n(z^-)^m
\leftrightarrow\cr \textrm{(or~Wick)} &\textrm{} 
&(\hat{z}^+)^n(\hat{z}^-)^m\cr &&\cr
\noalign{\hrule}
&&\cr
\textrm{Antinormal} &e^{-\frac{1}{2}\zeta ^+\zeta ^-} &(z^+)^n (z^-)^m
\leftrightarrow\cr (\textrm{or~Anti-Wick}) &\textrm{} 
&(\hat{z}^-)^m(z^+)^n\cr &&\cr
\noalign{\hrule}
&&\cr
\textrm{Symmetric} &\cos\frac{1}{4}[(\zeta ^+)^2-(\zeta ^-)^2] &q^np^m
\leftrightarrow\cr &&\frac{1}{2}(\hat{q}^n\hat{p}^m+\hat{p}^m\hat{q}^n)\cr 
&&\cr
\noalign{\hrule}
&&\cr
\textrm{Born-Jordan} &\frac{\sin\frac{1}{4} [(\zeta ^+)^2-(\zeta ^-)^2]}
{\frac{1}{4} [(\zeta ^+)^2- (\zeta ^-)^2]} &q^np^m\leftrightarrow\cr
&& \frac{1}{m+1} \sum_{j=0}^{m} 
\widehat{p}^{m-j}\widehat{q}^n\widehat{p}^j\cr &&\cr}
 \]

\newpage  
Suppose now that $\Omega (\zeta )\neq 0$  for all $\zeta \in \mathbb 
R^{2d}$. Then 
 $\tilde{f} (\zeta) /\Omega (\zeta )$ is meaningful for
  $f\in \mathcal S'(\mathbb R^{2d})$ 
 if and only if $1/\Omega$ 
is a multiplier in $\mathcal S(\mathbb R^{2d})$. In such a case $f^\Omega 
$ 
is called the {\it strict 
$\Omega $-symbol} of the distibution $f\in \mathcal S'(\mathbb R^{2d})$. 
E.g. every $f\in \mathcal S'(\mathbb R^{2d})$
 has strict standard , antistandard and normal symbols. More generally $f^
 {\Omega }$ is called the {\it strict $\Omega$ -symbol} of $f$
if only $\Omega (\zeta )\tilde{f} ^\Omega (\zeta ) \in \mathcal S' 
(\mathbb R 
^{2d}).$ \smallskip 

Of course not every Weyl operator has either antinormal, or symmetric, 
or Born-Jordan symbol.

\paragraph{Quasi-polynomials.}
Define (c.f.~\cite{shu}, Appendix 2) for 
$m=(m_1,m_2),\   r=(r_1,r_2),\  r_1 \geq 0,\  r_2<1/2$,
 the class $S(m,r)$ of 
{\it quasi-polynomial} $f=\{ (f_\hbar 
(z),\ 0<\hbar \leq \hbar (f) \}$ in 
$\mathcal S (\mathbb R ^{2d})$ such that 
\[
\partial _z^\alpha  f =\mathcal O _\alpha (1)(1+|z|)^{m_1 - r_1 |\alpha 
|}  
\hbar ^{m_2 - r_2 |\alpha |}  
\] 
wherein $\partial _z = \partial /\partial z$  and $\alpha $ is a
 multiindex. \newline 
As usual,
$S(-\infty ) = \cap S(m,r),\quad S(\infty ) = \cup  S(m,r)$. \smallskip

A quasi-polynomial $f$ is said to be  {\it asymptotic } to a series 
$\sum _{\alpha \geq \mu }{}  f^\alpha $ 
$$
f\simeq \sum_ {\alpha \geq \mu }{} f^\alpha 
$$
with $f^\alpha \in  S(m_\alpha,r_\alpha), \quad 
m_\alpha \searrow -\infty , \quad r_\alpha \searrow -\infty $ 
if for all $\nu$  
\[ 
f- \sum_{\alpha < \nu} f^\alpha \in  S(m_\nu,r_\nu).
\] 
The classical Borel-H\H{o}rmander construction leads to the following

\paragraph{Proposition.}
 {\it  For every $f\in S(m,r)$ and a formal $\Omega $ there is a $g\in 
S(m,r)$  
 asymptotic to $[1/\Omega (\frac{\hbar }{i} \partial _z )]f $.} \smallskip 
  
  Such  function $g$ is called an {\it asymptotic symbol} $f^{\Omega }$ of 
$f$. 
 It is defined $mod\  S(-\infty).$
 
 \paragraph{$\Omega$-products of quasi-polynomials.}
  
  If $f_j$ are quasi-polynomials then  $\hat{f_j} $ act from 
   $\mathcal G (\mathcal H)$ to $\mathcal G (\mathcal H)$ 
    and therefore from $\mathcal G' (\mathcal H)$ to $\mathcal G' 
    (\mathcal H)$, so that  $\hat{f} = \hat{f} _1 \hat{f} _2 \dots 
    \hat{f} _N$ is well defined.
    Actually $f$ is quasi-polynomial, and
  \[
  f^\Omega (z) = \int \mathcal K ^\Omega(z - z_1, \dots z - z_N) 
  \prod_{j=1}^{N} f_j^{\Omega} (z_j) d\lambda _{\hbar} z_j,
   \]
  wherein
   \[
  \tilde{\mathcal K} ^\Omega (\zeta _1, \dots ,\zeta _n)
    \simeq \frac {\prod_j \Omega (\zeta _j)}{\Omega (\sum _j \zeta _j)} 
\exp \left 
  \{ \frac{1}{2}  \sum_{j<k} [\zeta _j, \zeta _k] \right \}, 
   \]
   (=, if $\Omega$ is strict). Generally, the integral is distributional, 
   but  absolutely converges for some  $\Omega $, e.g., the  normal one.

  The integral representation entails the asymptotic expansion
   \[
   f^\Omega (z) \simeq \tilde { \mathcal K}^{\Omega }\left ( 
\frac{\hbar}{i} 
    \pder{}{z^+} ,- \frac{\hbar}{i} \pder{}{z^-} \right ) \prod_j f^\Omega 
(z_j) 
    | _{z_j =z},
    \]
   wherein
 $ \pder{}{z^+} = \frac{1}{\sqrt 2} (\pder{}{q} + \frac{1}{i} \pder{}{p} 
), \quad 
  \pder{}{z^-} = \frac{1}{\sqrt 2} (\pder{}{q} - \frac{1}{i} \pder{}{p}).$ 
  
  \paragraph{Trace.} A \emph{ density operator} $\hat{\rho}: \mathcal G' 
  \rightarrow \mathcal G $ is a Weyl operator with $\rho \in \mathcal 
  S (\mathbb R^{2d})$. The operator trace of $f^\Omega \rho ^\Omega $ is 
well defined and may be evaluated for strict $\Omega $ via the 
\emph{Trace formula}:
  Tr$ (\hat{f}^\Omega \hat{\rho } ^\Omega )=
  <f^\Omega |\rho ^\Omega> $ .

 \section {Main Theorem.}

A quasi-polynomial $f\in S(m,r),\  m>0,$ 
 is called {\it apt} if  
for sufficiently small $\hbar $
it satisfies the following three  conditions uniformly:
\begin{itemize}
  \item  \emph{Quasi-dissipativity}: $\emph{Re}(if) >\delta $, a constant .
  \item  \emph{Hypoellipticity}: for all multi-indices 
  $\alpha $ and $0\leq t',t'' \leq T$
   \[
   \partial _z^\alpha f(t'',z) =\mathcal O _\alpha (1) |if(t') - \delta |
   (1+|z|)^{ - r_1 |\alpha |}  \hbar ^{ - r_2 |\alpha |}.  
   \]
  \item  \emph{t-Continuity} of f(t,$\cdot $ ) in $S(m,r)$.
   \end{itemize} 
   {\it Law of Inertia}: If $f$ is apt then all its asymptotic symbols
  $f^\Omega $ are apt as well, albeit on different intervals of $\hbar $. 

  Also if $f$ is hypoelliptic and \emph{real} then $\hat{f}(t)$ are  
  essentially self-adoint  (c.f.~\cite{shu}, Proposition A2.1) in 
$\mathcal H$.

\paragraph{ Main Theorem.}
 {\it  If $f$ is an apt quasi-polynomial,  
  then for sufficiently small $\hbar$ 
\begin{description} 
 \item  (1) The Cauchy Problem 
  $$\fder{\psi}{t} + \hat{f} (t)\psi (z,t) = 0,\quad \psi (z,0) = 
\psi_0,\quad 
0\leq t\leq T, $$
is proper on $\mathcal H $ relative $\mathcal G (H)$. 
  \item (2) The evolution operator is the strong product integral
   \[
 U(t'',t') =\prod_{t''\geq t \geq t'} [ (\mathbf 1 + \frac{idt}{\hbar} {f} 
   (t,\cdot ))^{-1})]\  \widehat{\  }.
  \] 
  \item (3) A strict $\Omega$ -symbol $u^\Omega (t'',t',z)$ of the 
evolution 
   operator is  the limit in $\mathcal S' (\mathbb R ^{2d})$
   of 
    the strict $\Omega $-symbols $u^\mathcal P (t'',t',z)$ 
   of the partial operator products
      $\prod_{t''\geq t_j \geq t'} [ (\mathbf 1 + \frac{i\Delta 
t_j}{\hbar} {f} 
   (t_j,z))^{-1})]\  \widehat{\  }$ as $|\mathcal P| \rightarrow 0$. 
   
     \end{description}
   }

\paragraph{Proof (outline).} We apply the $\Omega $-calculus along with 
the theory 
of Abstract Cauchy Problems (c.f.~\cite{Fat}) and 
the theory of Finite Difference Methods for Initial Value Problems
(c.f.~\cite{Ric}) with the terminology thereof. 

{\it The following statements
 hold for various intervals of positive $\hbar$.} 
By the Law of Inertia, the \emph{anti-normal} symbol of $f$ is 
quasi-dissipative. Then 
 (c.f.~\cite{shu}, Proposition 24.1) the real part of
 \mbox{$<\! \psi |\delta _1 \mathbf 1 +i\hat{f} (t) |\psi \! >$}  is 
greater 
 than \mbox{$\gamma 
\! <\psi |\psi \!>$} with some constants $\gamma >0$ and $\delta _1$. It 
is 
  safe to  assume that  $\delta _1 =0$. Together with the hypoellipticity 
(c.f.~\cite{shu}, Theorem 25.4)  this entails that
   $|\!|[ \lambda  \mathbf 1 +i\hat{f} (t)]^{-1}| \!| <1/\lambda $ for 
positive
   $\lambda $ so that the operators
$\hat{f} (t)$ are a (1,0)-stable family in $\mathcal H$.

When both $\psi $ and $\hat{f} (t)\psi$ belong to $\mathcal H$ for some 
$t=t_0$ then it is so for all $t$ by the virtue of the  hypoellipticity.
 The space $\mathcal F$ of all such  $\psi $ is dense in $\mathcal H$ and 
 is a Hilbert space relative the new Hermitean 
 product $<\! \psi |\psi \!>_0=<\! \hat{f} (t_0)\psi |\hat{f} (t_0)\psi 
\!>$. Now the 
 $\hat{f} (t): \mathcal F \rightarrow \mathcal H$ form a $t$-continuous 
 family of bounded operators. Moreover, \mbox{$<\!\hat{f}(t)\psi |\psi 
\!>_0$}
 =\mbox{$<\! \hat{g}(t)\psi |\psi \! >$} with 
 $\hat{g}(t)$=\mbox{$\hat{f}(t_0)^\dagger  \hat{f}(t_0)\hat{f}(t)$} 
 so that $g$ is apt again and thus (as above) is (1,0)-stable in  
$\mathcal F$.
 By the Hille-Yosida theorem~\cite{hil}, $\hat{f} (t)$ generates 
 for every $t$ a contractive operator semi-group in $\mathcal F$.
  Since the family \mbox{$\hat{f} (t)-\hat{f} (T) $} has similar 
properties,
   the theorem 7.7.13
of~\cite{Fat} establishes (1), the \emph{properness} of the Cauchy 
problem. 
 
This leads to a preliminary  Product Integral representation
 $U(t'',t') =\prod_{t''\geq t \geq t'} [\mathbf 1 +\frac{idt}{\hbar} 
\hat{f}^{-1}]$
 (c.f. the proof of the theorem 7.7.5 of~\cite{Fat}). It implies 
 the Product Integral representation (2)
 via the Lax Equivalence Theorem~\cite{Ric} whereby the required 
consistency 
 is checked via the \emph{Weyl}  calculus.

The last statement (3) follows from the \emph{trace formula}.

 \section{Path Integrals in Quantum Field Theory.}
 \paragraph{Infinite dimensional phase spaces.}
 In the case of $d=\infty$ there are non-isomorphic phase spaces and  
  the symplectic structures usually appear with extra features.
 
 Our phase  space is based on a separable Frechet nuclear space $\mathcal 
Z$ 
 over $\mathbb C$ with a ``dotless'' hermitian product $zw^*$. 
 
 If $\mathcal H$ is the corresponding
 Hilbert space completion of the $\mathcal Z$,  
 and $\mathcal Z^*$ is the corresponding anti-dual of $\mathcal Z$ then
 $$\mathcal Z \hookrightarrow \mathcal H \hookrightarrow \mathcal Z^*$$ 
 is a Gelfand nuclear triplet.
 
 The phase space is $\mathcal Z$ taken over $\mathbb R$ with the 
 symplectic form $-\mathrm{Im}(zw^*)$.
 It is also a pre-Hilbert space with the scalar product 
$\mathrm{Re}(zw^*)$.
 
 \paragraph{Complex Gaussian rigging.} (C.f.~\cite{Hid}.)
 \emph{The Gaussian measure} $\gamma _\hbar $ of covariance $1/\hbar$ is 
defined 
 via its characteristic function
 $$
 \int_{\mathcal Z^*} e^{i\mathrm{Re} (z\zeta ^*)} d\gamma _\hbar (\zeta 
 ^*) = e^{-zz^*/2\hbar },\quad z\in \mathcal Z,
 $$
 so it stands for the  non-existent $(\frac{\hbar}{2\pi} )^{\infty }
  \exp (-\hbar \zeta ^* \zeta ^* ) d\zeta ^* $.
 
 The Bargmann-Segal space $[\mathcal H]$ is the closure of the subspace of 
 the continuous 
 complex analytic polynomials on $\mathcal Z ^*$ in $\mathcal
  L ^2 (\mathcal Z ^*)$. Its elements are entire functions $h(z^*)$ of 
  order 2 and type $<\hbar /2$:
  $$
  h(z^*)=\mathcal O (e^{ \hbar p(z^*)^2 /2})
  $$
  for some dual semi-norm $p$ on $\mathcal Z^*$.

  Let $[\mathcal Z]$ denote the space of entire functions $h(z^*)$ of 
  order 2 and minimal type, and $[\mathcal Z]^*$ 
  denote the space of entire functions $h(z^*)$ of order 2 and maximal 
type.
  Then $[\mathcal Z]$ is naturally a separable nuclear Frechet space, and
  $[\mathcal Z] ^*$ its anti-dual. 
  Thus
 \[
 [\mathcal Z] \hookrightarrow [\mathcal H] \hookrightarrow [\mathcal Z] ^*
 \]
 is another Gelfand triplet, the Complex Gaussian rigging  of the triplet
 \[
  \mathcal Z \hookrightarrow \mathcal H \hookrightarrow \mathcal Z ^*.
 \]
 The \emph{coherent states} $e_w (z^*):=\exp (-\hbar wz^*), \ w\in 
 \mathcal Z$, form a total (overcomplete) set in $[\mathcal Z]$.
  
\paragraph{Free bosonic field over $\mathcal Z$ in $[\mathcal Z]^*$.} 
  Let $z\rightarrow \bar{z} $ be an antilinear conjugation on 
  $\mathcal Z$ and correspondingly on $\mathcal Z^*$. Set

\[
 z^+=z/\sqrt 2 \in \mathcal Z,\ z^- = \bar{z}/\sqrt 2 \in \bar{\mathcal Z}.
 \]

 The operators $\hat{z} ^+$ and $\hat{z} ^-$ are defined on $f\in 
[\mathcal 
 Z ]^*$ as
 \[
 \hat{z}^+ f(\zeta ^*) = (z\zeta ^*)f(\zeta ^*), \ \hat{z} ^- =\hbar 
 \partial _z f(\zeta ^*).
 \]
  They represent the Canonical Commutation Relations (CCR):
 \[
 [\hat{z} ^-,\hat{z} ^+]=\hbar \mathbf 1.
  \]
The coherent states are entire vectors for the CCR.  
\paragraph{Wick Operators} are, by definition, the continuous linear 
operators $W$
from $[\mathcal G]$ to $[\mathcal G]^*$.
 The \emph{Wick symbol} of $W$ is 
 \[
 w(z^+,z^-):=e^{-z^+z^-}\int _{\mathcal Z^*}[We_{z^+}(\zeta ^-)]e_{z^-}    
(\zeta ^+)d\gamma _{\hbar }(\zeta ),\quad z\in \mathcal Z.  
\]

The Wick symbol $w$ is an entire function on $\mathcal Z \times \mathcal 
Z$,
 so that the operator $\hat{w}: = w(-\hat{z}^+,\hat{z}^-)$ is 
well defined on the coherent states and $W=\hat{w} $ and $W$ is its 
closure.  
 
\paragraph{$\Omega$ -symbols.}
Consider a formal complex power series $1+\sum_{|\alpha |>0}z^{\alpha }$ 
on 
 $\mathcal Z^*$.
The formal $\Omega $-symbol of $w(\zeta )$ is  
\[
w^{\Omega }:= \left [1/\Omega (\frac{\hbar }{i}\partial _z)\right ]w(\zeta 
).
\]

\paragraph{Quasi-polynomial}$w(\zeta )$ is the family $\{ w_{\hbar } (\zeta
) :0<\hbar \leq \hbar (w)\}$ such that for a dual semi-norm $p$ on 
$\mathcal
 Z ^*$
 \[
 \partial _z ^{\alpha }w(\zeta ) = \mathcal O _{\alpha ,p}(1)(1+p(\zeta ))
 ^{m_1 -r_1 |\alpha |} \hbar ^{m_2 -r_2|\alpha |}.
\]
The class of such families is denoted $S(m,r),\ m=(m_1,m_2),r=(r_1,r_2).$
\paragraph{Weyl symbols and operators.}

The Weyl symbols $f(z)$ of $w(z)$ correspond to 
$\Omega (z) = \exp (-\frac{1}{2} z^+z^-)$ (c.f.the table above):
\[
f(z) = \exp \left [\frac {\hbar ^2}{2} \frac {\partial _2}
{\partial z^+ \partial z^-}  \right ] w(z). 
\]
Conversely
\[
w(z) =\int_{\mathcal Z ^*}f(z-\zeta )d\gamma _{\hbar ^2 /2} (\zeta ).
\]
(Note: not every Wick operator has a strict Weyl symbol.)\smallskip 

The corresponding Wick operators are  the \emph{Weyl operators} 
$\hat{f}$.
The Borel-H\H{o}rmander constructions for a countable fundamental family 
of dual 
gaussian semi-norms followed by the Cantor diagonal trick imply that for 
every $\Omega $ and Weyl $f\in \mathcal S (m,r)$ there is $g\in \mathcal S 
(m,r)$ asymptotic to $f^{\Omega }$.\smallskip

 The $\Omega $-symbols of the operator product $\hat{f}_1 \hat{f}_2 ... 
 \hat{f}_N$ have  
 the asymptotic expansions just as in the case $d<\infty $. 
   However their integral representations are known rarely.
   Fortunately, for the \emph{normal symbols} $w$ of $\hat{w}_1\hat{w}_2 
...
   \hat{w}_N
   $
   \[
   w(z) = \int \prod_{j=1}^{N} e^{z_j^- z_{j-1}^+} w_j (z_j^-,z_{j-1}^+) 
    \prod_{j=1}^{N-1} d\gamma_{\hbar} (z_j^-,z_{j-1}^+),\quad 
    z_0^+:=z^+,\ z_N^-:=z^-.
    \]  
  
Finally, as in the finite-dimensional case, \emph {the Main Theorem} holds
in infinite dimensions (with the same proof) 
at least for the strict Wick symbols of evolution operators in 
$G=[\mathcal H]$. In the latter case
the symbol approximations are absolutely convergent multiple integrals
 with respect to 
$d\gamma _\hbar $ over the infinite-dimensional phase space $\mathcal Z$. 
\smallskip

      \section{Conclusion and outlook.}
   \begin{enumerate}
        
  \item The  phase space Path Integral (according to L.Shulman, ''a 
 difficult form'' of the path integral) 
was originated in different ways by Feynman himself ~\cite{Fey} in 1951 
  and by Tobocman ~\cite{Tob} in 1956. The Coherent State discretization 
  was introduced in 1960 by Klauder~\cite{kla} in the Schr\H{o}dinger 
representation
   and in 1962 by Schweber~\cite{sch} in tne Bargmann-Fock 
represenatation. In the 70's 
 Berezin~\cite{Ber} considered various discretizations 
  on the basis of pseudo-differential analysis. However no convergence of 
  the discretizations  has been proved until now.
  
   On the other hand
  Daubeshies \& Klauder~\cite{Dau} have established in 1984 that  a wide 
class of 
 coherent state path integrals (essentially with self-adjoint polynomial 
 hamiltonians) on a flat finite-dimensional phase space  may are limits
 of  Wiener Integrals on the space of paths in the phase space. They even 
 suggested that  Feynman type time-slicing construction is impossible 
 for the phase space Path Integrals. 

 \item  We have presented a \emph{rigorous time-slicing} Phase Space Path 
Integral
     construction for the symbols of the Evolution Operators  with wide 
variety 
     of smooth hamiltonians both in finite and infinite degrees of 
freedom.    
 The  convergence is established only for small $\hbar $, 
   in agreement with  the postulated semi-classical nature of the Path 
Integral
    which relates the classical and quantum dynamics.
  \item  According to the $\Omega $-calculus, the discretizations 
  of the Path Integral are distributional multiple integrals.
   E.g., as mentioned in the Introduction, the traditional discretization 
   of the Phase Space  Integral comes from the standard $\Omega 
$-calculus. 
   Similarly, 
  the Coherent State Path Integral discretization is associated with the
   normal $\Omega $-symbol in which case the multiple integrals are 
   absolutely convergent.
 \item  The last statement in the Main Theorem is equivalent to a modified 
   DTF-ansatz: the $\Omega $-symbol of the short time propagator is 
   approximately equal to $[\mathbf 1 +\frac{i}{\hbar}f(t,z)\Delta 
t]^{-1}$. 
   However, in the case of the normal $\Omega$-calculus (because of the 
absolute 
   convergence) one may consistently replace it with the more customary 
ansatz 
 $\exp [-\frac{i}{\hbar} f(t,q,p)\Delta t]$.
       
  \item Our Path Integrals are ``pathless'', in agreement with the 
    Uncertainty Principle: no quantum path in the phase space 
    (c.f.~\cite{Kla} for an illuminating somewhat different point of 
view). 
  Yet they are semiclassical in the following sense: the principal terms 
in 
  the $\hbar $-expansions 
  of the partial products symbols are the backward Euler approximations 
  of the  corresponding  classical  Hamilton-Jacobi equations.
   \item  We have rigorized  the $\Omega $-calculus of 
    Agarval \& Wolf~\cite{Aga} to justify numerous  Path Integral 
    discretizations and as an important techniques.
    However in the  infinite degrees we have been able to prove the 
     convergence only for normal (Wick) symbols. 
     
     Actually the formal $\Omega $-calculus is a special case of the 
formal 
    $\ast $-calculus~\cite{Bay}. Since on the finite-dimensional flat 
     symplectic space all $\ast $-products are formally 
    equivalent, our results yield a {\it construction of the formal
     \mbox{$\ast $-exponential}},
    a solution of a well known problem (c.f~\cite{sha} for an 
    interpretation of the Evolution Operator symbol as a 
$\ast$-exponential).

  \item  Most of the other mathematical  interpretations of the Path 
 Integral are primarily in terms of various 
   distributional integrals on the paths in the configuration  space : 
first, 
 by Kac~\cite{Kac} via analytic continuation to a Wiener integral (the  
Feynman-Kac formula), followed by  DeWitt-Morette~\cite{Mor}  
in terms of prodistributions,  
by  Albeverio and H\o eg-Krohn~\cite{Alb}  in terms of the Parseval 
equation 
for the oscillatory Gaussian integrals, and by 
 Hida \& Streit~\cite{Hid} in terms of White Noise disributions.
 Notably, these Path Integrals are associated  
only with the Schr\H{o}dinger hamiltonians (essentially)
  of quadratic growth, with the presumed ``Feynman measure'' built
 from the kinetic energy term.
\end{enumerate}

  \end{document}